\documentclass[12pt,oneside]{amsart}

\usepackage{amsthm}
\usepackage{mathrsfs}
\usepackage{enumerate}
\usepackage{amsfonts}
\usepackage{verbatim}
\usepackage{amsbsy}
\usepackage{amsmath}
\usepackage{amssymb}
\usepackage{MnSymbol}

\newtheorem{theorem}{Theorem}[section]
\newtheorem{lemma}[theorem]{Lemma}
\newtheorem{proposition}[theorem]{Proposition}
\newtheorem{conjecture}[theorem]{Conjecture}
\newtheorem*{burning-conjecture}{Burning number conjecture}
\newtheorem*{claim}{Claim}
\newtheorem{question}[theorem]{Question}

\theoremstyle{definition}
\newtheorem{definition}[theorem]{Definition}
\newtheorem*{notation}{Notation}
\theoremstyle{remark}
\newtheorem{remark}[theorem]{Remark}

\newtheorem{case}{Case}

\numberwithin{subcase}{case}
\numberwithin{subsubcase}{subcase}

\DeclareMathOperator{\summation}{sum}

\makeatletter    %Need to allow MSC year to be shown as 2020
\@namedef{subjclassname@2020}{%
	\textup{2020} Mathematics Subject Classification}
\makeatother

\subjclass[2020]{05C85, 68R10, 91D30}

\begin{document}

\title{Burnability of Double Spiders and Path Forests}

\author{Ta Sheng Tan}
\address{Institute of Mathematical Sciences, Faculty of Science, Universiti Malaya, 50603 Kuala Lumpur, Malaysia}
\email{tstan@um.edu.my}

\author{Wen Chean Teh\!\! $^*$}
\address{School of Mathematical Sciences, Universiti Sains Malaysia, 11800 USM,\linebreak Malaysia}
\email{dasmenteh@usm.my}

\keywords{Spread of social contagion; Burning number conjecture; Graph algorithm; Double spider; Path forest}

\begin{abstract}
 The burning number of a graph can be used to measure the spreading speed of contagion in a network.
 The burning number conjecture is arguably the main unresolved conjecture related to this graph parameter, which can be settled by showing that every tree of order $m^2$ has burning number at most $m$.
 This is known to hold for many classes of trees, including spiders - trees with exactly one vertex of degree greater than two.
 In fact, it has been verified that certain spiders of order slightly larger than $m^2$ also have burning numbers at most $m$, a result that has then been conjectured to be true for all trees.
 The first focus of this paper is to verify this slightly stronger conjecture for double spiders - trees with two vertices of degrees at least three and they are adjacent.
 Our other focus concerns the burning numbers of path forests, a class of graphs in which their burning numbers are naturally related to that of spiders and double spiders.
 Here, our main result shows that a path forest of order $m^2$ with a sufficiently long shortest path has burning number exactly $m$, the smallest possible for any path forest of the same order.
\end{abstract}

\let\thefootnote\relax\footnotetext{$^*$\! \!Corresponding author}

\maketitle

\section{Introduction}

Graph burning is a discrete-time process introduced by Bonato, Janssen, and Roshanbin~\cite{bonato2014burning, bonato2016how, roshanbin2016burning} that can be viewed as a simplified model for the spread of contagion in a network.
Given an undirected finite graph $G$ without loops and multiple edges, each vertex of the graph is either \emph{burned} or \emph{unburned} throughout the process.
At first, every vertex of $G$ has the unburned status. Sequentially, 
a \emph{burning source} (or simply \emph{source}) is placed at an unburned vertex to burn it
at the start of every round $t\ge 1$.
If a vertex is burned in round $t-1$, then in round $t$, each of its unburned neighbors becomes burned.
A burned vertex is assumed to remain burned throughout the burning process.
The process terminates when all vertices of $G$ have acquired the burned status, in which case we say the graph $G$ is \emph{burned}.
 The least number of rounds required to complete the burning process is called the \emph{burning number of $G$} and it is denoted by $b(G)$.

While the burning numbers of graphs do not satisfy general hereditary property in the sense of taking subgraphs, Bonato, Janssen, and Roshanbin~\cite{bonato2016how} showed that in order to burn a connected graph, it suffices to burn its spanning trees.
Therefore, much focus has been spent on trees in the study of graph burning.
As one of their initial results, they also determined the exact burning numbers of paths (and thus cycles or hamiltonian graphs), showing that $b(P_m) = \lceil\sqrt{m}\rceil$ with $P_m$ being the path of order $m$.
This marked the beginning of their currently still unresolved main conjecture on graph burning.

\begin{burning-conjecture}[\cite{bonato2016how}]
 If $G$ is a connected graph of order $m$, then $b(G)$ is at most $\lceil\sqrt{m}\rceil$. 
\end{burning-conjecture}

In the literature on graph burning, a graph that satisfies the burning number conjecture is said to be \emph{well-burnable}.
So paths and hamiltonian graphs are well-burnable.
As remarked above, the burning number conjecture holds if all trees are well-burnable.
Classes of trees known to be well-burnable include spiders \cite{bonato2019bounds, das2018burning} and caterpillars \cite{hiller2019burning, liu2020burning}.
Here, a \emph{spider} is a tree with one vertex of degree greater than two, and a \emph{caterpillar} is a tree where a path remains after deleting all vertices of degree one.
For general connected graphs $G$ of order $m$, since the initial bound of $b(G)\le 2\sqrt{m}-1$ in \cite{bonato2016how}, some attempts~\cite{bessy2018bounds, land2016upper} have been made towards proving the burning number conjecture, with the currently best known upper bound being roughly $\frac{\sqrt{6}}{2}\sqrt{m}$ by Land and Lu~\cite{land2016upper}.

Beyond being well-burnable, there are of course trees of order greater than $m^2$ having burning numbers at most $m$, with the simplest such graphs being the stars, where every one of them has burning number exactly two.
Together with the intuition that a tree that deviates from a path should be easier to burn, the following slightly stronger conjecture was made in \cite{tan2020graph}.
For convenience, if the burning number of a graph is at most $m$, then
we say the graph is \emph{$m$-burnable}. 

\begin{conjecture}\label{Tree Conjecture}
 Let $m>n\ge 2$. If $T$ is a tree with $n$ leaves and its order is at most $m^2+n-2$, then $T$ is $m$-burnable.
\end{conjecture}

This conjecture obviously holds for stars and paths, and it has been verified in the same paper for spiders.
Calling a spider with $n$ leaves an \emph{$n$-spider}, we note that for any $m,n\ge 2$, there are $n$-spiders of order $m^2+n-1$ that are not $m$-burnable, which implies that the bound in the conjecture is tight.
While not all $n$-spider of order $m^2+n-2$ are $m$-burnable for the case when $m \le n$, those that are not $m$-burnable must contain, as a subgraph, the $m$-spider of order $m^2+1$ such that the distance between any two leaves is exactly $2m$. 
These are the results on the burnability of spiders in \cite{tan2020graph}.

One of the purposes of this paper is to further verify Conjecture~\ref{Tree Conjecture} for the next natural class of trees called double spiders - the union of two spiders with an edge joining their respective vertices of maximum degree.
More precisely, a \emph{double spider} is a tree that has two special adjacent vertices called \emph{heads} with the property that every other vertex has degree at most two. 
Note that in a double spider, every leaf that is not a head must be connected to the closest head by a unique path, which we will call an \emph{arm} of the double spider.
An \emph{$n$-double spider} is a double spider with $n$ arms.
With this definition, paths and spiders can be viewed as double spiders.
We also remark that while most $n$-double spiders have exactly $n$ leaves, those where all the arms are joined to the same head have $n+1$ leaves.
%This flexibility is helpful in proving our results on the burnability of double spiders.
Our main result here verifies Conjecture~\ref{Tree Conjecture} for double spiders.

\begin{theorem}\label{double_spiders}
 Let $m>n\ge 2$. If $T$ is an $n$-double spider of order at most $m^2 + n - 2$, then $T$ is $m$-burnable.
\end{theorem}

For the case when $m\le n$, the same examples as those for spiders show that not all $n$-double spiders of order $m^2 + n - 2$ are $m$-burnable. 
To complete our study on the burnability of double spiders, we show in Theorem~\ref{double_m<n} that other than similar such cases, every $n$-double spider of order $m^2 + n -2$ is $m$-burnable.
We note that the bound of $m^2+n-2$ is again tight for any $m,n\ge 2$, simply by attaching $n-2$ leaves to two adjacent vertices on a path of order $m^2 + 1$.

A \emph{path forest} is a disjoint union of paths.
The study of graph burning of spiders is naturally related to that of path forests, as observed in the previous studies~\cite{bonato2019bounds,das2018burning,tan2020graph} on burning spiders.
In particular, it was shown in \cite{tan2020graph} that while every path forest consisting of $n$ paths and of order $m^2 - (n-1)^2$ is \mbox{$m$-burnable}, the only path forest with $n$ paths and of order 
one larger than $m^2 - (n-1)^2$ that is not \mbox{$m$-burnable} contains $n-1$ independent edges as components.
This prompts us into investigating the burnability of larger path forests.
Could it be possible that the only larger path forests that are not $m$-burnable is a small collection of nicely described path forests?

Therefore, our second focus in this paper is on the burnability of path forests.
To this end, we first slightly refine the above result in \cite{tan2020graph}, showing that if $T$ is a path forest with $n$ paths such that its shortest path has order at least six and $$|T|\le m^2 - (n-1)(n-2)+1, $$ then $T$ is $m$-burnable.
(See Theorem~\ref{n-paths} for the precise version of the result.)

Of course, no path forest of order larger than $m^2$ is $m$-burnable. 
We observe from existing results that a path forest consisting of two paths and of order $m^2$ is \mbox{$m$-burnable}, provided that the shorter path is not just an edge.
What about path forests with more paths?
Is it true that every path forest of order $m^2$ with sufficiently long paths is $m$-burnable?
While we are not able to quantify ``sufficiently long'' here, even if we fix the number of paths in the path forest, we show that the answer to this question is affirmative.

\begin{theorem}\label{long-path-forests}
 Let $n\geq 2$. There exists $L\in \mathbb{N}$ with the following property.
For every path forest $T$ with $n$ paths, if the shortest path of $T$ has length at least $L$ and $|T| = m^2$, then $T$ is $m$-burnable.	
\end{theorem}

As one may notice, the above problem on burning path forests is equivalent to the following partition problem.
Given $n$ positive integers $l_1,l_2,\ldots,l_n$ summing to $m^2$, we wish to decide if the set of the first $m$ odd positive integers can be partitioned into $n$ sets $S_1, S_2, \ldots, S_n$ such that for every $1\le i \le n$, the sum of the numbers in $S_i$ is equal to $l_i$.

The plan of the paper is as follows.
In Section~\ref{path-forests}, we observe some simple results on path forests that are useful for burning double spiders, and we investigate \mbox{$m$-burnable} path forests with $n$ paths and of order larger than $m^2-(n-1)^2 + 1$.
Then in Section~\ref{double}, we prove our results on burning double spiders.
Our result on the burnability of path forests with sufficiently long shortest paths will be the content of Section~\ref{path-forests-long}.
Finally, we mention some remarks and open problems in Section~\ref{conclusion}.

While our work focuses on burning trees and path forests, we remark that graph burning of various other classes of graphs has been studied, such as graph products~\cite{mitsche2018burning}, hypercubes~\cite{mitsche2018burning}, and random graphs~\cite{mitsche2017burning}.
Since the introduction of graph burning less than a decade ago, a considerable amount of work on this topic has produced many results and algorithms, including those on variants of it.
An excellent survey on the topic of graph burning can be found in \cite{MR4233796}.

We recall some terminology in the context of graph burning.
In a burning process, a \emph{burning sequence} is the sequence of vertices at which the burning sources are placed  in each round. The shortest such sequence is said to be  \emph{optimal}.
So the length of an optimal burning sequence of a graph is the burning number of the graph.
For a path forest, its \emph{path orders} indicate the respective order of each of its paths.
For an arm of a double spider, we usually consider its vertices in order, with the vertex next to the head it is joined to as the first vertex and so forth while the $0$th vertex is the head to which the arm is joined.

% This is the edited sections on path forests
%
%
%
%
%

\section{Path Forests}\label{path-forests}

In this section, we study the burnability of path forests.
We begin with the burnability of some special classes of path forests, which will be useful in proving the burnability of double spiders in Section~\ref{double}.
For the main result of this section, we obtain a strengthening of the following result in \cite{tan2020graph} to accommodate larger path forests but with more exceptional cases.

\begin{theorem}\cite{tan2020graph}\label{092021a}
Suppose $m\geq n\geq 2$ and let $T$ be a path forest with $n$ paths. If
$$\vert T\vert \leq m^2-(n-1)^2+1,$$
then $T$ is $m$-burnable unless equality holds and the path orders of $T$ are 
$$m^2-n^2+2, \underbrace{2,2,\dotsc, 2}_{n-1 \text{ times}}.$$
\end{theorem}

We first recall a simple lemma from \cite{tan2020graph} on the burnability of path forests with $n$ paths and of order at most $3n-2$.

\begin{proposition}[\cite{tan2020graph}]\label{n-paths-3n}
Suppose $n\geq 2$ and let $T$ be a path forest with $n$ paths.
 If $\vert T\vert\leq 3n-2$ and the shortest path of $T$ has one vertex, then $T$ is $n$-burnable.
\end{proposition}

This simple result was verified by a straightforward induction on the number of paths of $T$.
Using almost identical arguments by induction, we have the following simple burnability results for path forests.

\begin{proposition}\label{n-paths-linear}
Suppose $n\ge 2$ and suppose $T$ is a path forest with path orders $l_1\ge l_2\ge \cdots \ge l_n$.
 \begin{enumerate}[(i)]
  \item\label{n-paths-4n} If $\vert T\vert \leq 4n-4$ with $l_n = 1$ and $l_{n-1}\ge 2$, then $T$ is $n$-burnable.
  \item\label{n-paths-5na} If $ \vert T\vert\leq 5n-6$ with $l_n = 1$ and $l_{n-1} = 3$, then $T$ is $n$-burnable.
  \item\label{n-paths-5nb} If $\vert T\vert\leq 5n-1$ with $l_n \ge 3$, then $T$ is $(n+1)$-burnable.
 \end{enumerate}
\end{proposition}

\begin{proof} 
 For each of the three statements, the base case $n=2$ is straightforward.
 We suppose now that $n>2$ and the statements are all true for $n-1$.

 Let $T$ be a path forest as in \eqref{n-paths-4n}.
 It is obvious that $T$ is $n$-burnable when $l_1\le 3$, and so we consider the case when $l_1\ge 4$.
 Noting that $l_1\le 2n-1$ and $l_2+l_3+\cdots +l_n \le 4(n-1) - 4$, we see that $T$ is $n$-burnable since by removing the first path of $T$, the resulting path forest is $(n-1)$-burnable by induction hypothesis.

 For a path forest $T$ as in \eqref{n-paths-5na}, we must also have $l_1\le 2n-1$. So if $l_1\ge 5$, we similarly deduce that $T$ is $n$-burnable by induction hypothesis. Since $n\ge 3$, it is also obvious that $T$ is $n$-burnable when $l_1\le 4$.

 Finally, suppose $T$ is a path forest as in \eqref{n-paths-5nb}.
 Again, if $5\le l_1\le 2n+1$, $T$ is $(n+1)$-burnable by induction hypothesis, and if $l_1\le 4$, $T$ is also $(n+1)$-burnable where one of the paths can be burned using the last two burning sources.
 As $l_n\ge 3$, the only case left is when $l_1 = 2n+2$ and $l_2=l_3=\cdots=l_n=3$, which is clearly $(n+1)$-burnable where the first path can be burned using the first and the last burning sources.
\end{proof}

\begin{remark}
Suppose $T$ is a  path forest with path orders $l_1, l_2, \ldots, l_n$ and $T'$ is a path forest with path orders $l_1', l_2', \ldots, l_n'$ such that $l_i'\le l_i$ for each $1\le i\le n$. If $T$ is $m$-burnable, then $T'$ is also clearly $m$-burnable.
\end{remark}

Before we prove the main result of this section, we shall introduce some notation on path forests.
For $m\ge n\ge 2$, let $\mathcal{T}_{n,m}$ be the set of all path forests with $n$ paths where the path orders $l_1,l_2,\ldots,l_n$ are such that
\begin{enumerate}[(i)]
 \item $l_1 = m^2-(n-1)^2+1$ and $l_2=l_3=\cdots=l_n=1$, or
 \item $l_1= m^2-n^2+2$ and $2\le l_2, l_3,\ldots, l_n \le 3$, or
 \item $l_1= m^2-(n-1)(n+3)+1$ and $l_2=l_3=\dotsb = l_n=5$.
\end{enumerate}
(When $m=n$, only (i) is applicable.)

It is straightforward to verify that each of the path forests in $\mathcal{T}_{n,m}$ is not \mbox{$m$-burnable}.
For convenience, for a path forest $T$, we define $t_T$ to be the number of paths of order two in $T$ if the shortest path of $T$ has order two, and $t_T = 0$ otherwise.
So every path forest $T$ in $\mathcal{T}_{n,m}$ has order exactly $m^2-(n-1)(n-2)+1-t_T$.

The improvement over Theorem~\ref{092021a} states that if $T$ is a path forests with $n$ paths and the order of $T$ is at most $m^2-(n-1)(n-2)+1-t_T$, then $T$ is \mbox{$m$-burnable}, provided $T$ is not in $\mathcal{T}_{n,m}$.
The following lemma on path forests with three paths deals with the base case of this main result.

\begin{lemma}\label{3-paths}
 Let $m\ge 3$ and let $T$ be a path forest consisting of three paths.
 If $$|T|\le m^2-1-t_T\quad\mbox{and}\quad T\notin\mathcal{T}_{3,m},$$
 then $T$ is $m$-burnable.
\end{lemma}

\begin{proof}
 We prove by mathematical induction on $m$.
 Note that for the base case $m=3$, we must have $t_T=0$.
 Since $T\notin \mathcal{T}_{3,3}$ is a path forest of order at most eight, its path orders $l_1\ge l_2\ge l_3$ must satisfy $l_1\leq 5$, $l_2\leq 3$, and $l_3=1$, and thus $T$ is $3$-burnable.
 
 Now, let $m\geq 4$ and suppose the lemma is true for $m-1$.
 Let $T$ be a path forest with path orders $l_1\geq l_2\geq l_3$ such that $|T|\le m^2-1-t_T$ and $T\notin\mathcal{T}_{3,m}$.
 If $|T|\le m^2 - 3$, we see that $T$ is $m$-burnable by Theorem~\ref{092021a}. 
 So we will only need to consider the case when $|T|\ge m^2 - 2$, and in particular, we may assume $t_T\le 1$.
 
 \begin{case}
  $t_T=1$
 \end{case}
 In this case, $l_1+l_2+l_3 = m^2-2$, $l_3=2$, and $l_2\ge 4$.
 If $l_1\ge 2m+3$, the path forest $T'$ with path orders $l_1-(2m-1)\ge 4, l_2,l_3$ is $(m-1)$-burnable by induction hypothesis as $|T'|=(m-1)^2 - 2 = (m-1)^2 - 1 - t_{T'}$ and $T'$ is clearly not in $\mathcal{T}_{3,m-1}$.
 
 %If $l_1\le 2m + 2 \le (2m-1) + 5$, the first path of $T$ can be burned with the first and the $(m-2)$th burning sources, while the $(m-1)$th burning source can be used to burn the last path.
 %Now the remaining burning sources can burn a path of order $(m-1)^2 - 3 - 5$, and this is greater than $2m+2\ge l_2$, provided $m\ge 6$.
 
 Observe that we can have $l_1\le 2m + 2$ here only when $m\le 5$ as $m^2 - 2 > 2(2m+2) + 2$ whenever $m\ge 6$.
 %The only cases remain are when $m=4$ and $m=5$ with $l_1\le 2m+2$.
 When $m=4$, since $|T|=14$, the path orders of $T$ can be $(8,4,2)$, $(7,5,2)$, or $(6,6,2)$, each of which is clearly $4$-burnable. Similarly for $m=5$, the path orders of $T$ can only be $(12,9,2)$ or $(11,10,2)$, and so $T$ is $5$-burnable.
 
 \begin{case}
  $t_T=0$.
 \end{case}
 In this case, $m^2-2\le l_1+l_2+l_3\le m^2-1$ and $l_3\neq 2$. If $l_3 = 1$, the path forest $T'$ with path orders $l_1+1$ and $l_2+1$ has order at most $m^2$. Since $T\notin\mathcal{T}_{3,m}$, we have $l_2 \ge 2$, and so $T'$ is $m$-burnable by Theorem~\ref{092021a}.
 We can then deduce that the first two paths of $T$ can be burned using the first $m-1$ burning sources, and the last path can be burned using the last burning source.
 So we may assume that $l_3\ge 3$.
 
 Observe that $l_1\le 2m-3$ is possible only when $m=4$ as $m^2-2 > 3(2m-3)$ whenever $m\ge 5$.
 If $m = 4$, the only possible $T$ with $l_1\le 2m-3 = 5$ has path orders $(5,5,4)$, which is clearly $4$-burnable.
 So we may further assume that $l_1\ge 2m - 2$.
 
 Consider the path forest $T'$ obtained from $T$ as follows: delete a path of order $2m-2$ or $2m-1$ if there is such a path in $T$, or else delete $2m-1$ vertices from the first path of $T$.
 %Consider the path forest $T'$ obtained from $T$ by deleting $\max\{l_1, 2m-1\}$ vertices from its first path.
 Clearly,  $T$ is $m$-burnable  if $T'$ is $(m-1)$-burnable.
 Observe that either $T'$ is a path forest with two paths and  $\vert T'\vert \leq (m-1)^2$, or $T'$ is a path forest with three paths and $\vert T'\vert \leq (m-1)^2 - 1$.
 For the former, $T'$ is $(m-1)$-burnable by Theorem~\ref{092021a}, and for the latter, we see that $T'$ is $(m-1)$-burnable from induction hypothesis unless $T'\in\mathcal{T}_{3,m-1}$ or $T'$ has order $(m-1)^2 - 1$ with $t_{T'} = 1$.
 
 Since $l_3\ge 3$ and $T\notin\mathcal{T}_{3,m}$, for $T'$ to be a path forest in $\mathcal{T}_{3,m-1}$, it must either be the case that $l_1\in \{2m+1,2m+2\}, l_2\ge 4, l_3=3$ or $l_1 = 2m+4, l_2\ge 6, l_3=5$.
 It is now straightforward to verify that $T$ is $m$-burnable in either of the cases.
 In the latter case, we must have $m\ge 6$ since $m^2 - 1 \ge l_1+l_2+l_3\ge 2m+15$, and so the first path of $T$ can be burned with the second and the $(m-3)$th burning sources (which would burn $2m-3 + 7 = 2m+4$ vertices) and the last path of $T$ can be burned with the $(m-2)$th burning source, while the remaining burning sources are enough to burn the second path of $T$.
 The other case can be verified similarly. 

 The final case left is when $T'\notin\mathcal{T}_{3,m-1}$ is of order $(m-1)^2-1$ with $t_{T'} = 1$.
 We must then have $l_1 = 2m+1$ and $l_3\ge 4$, which is only possible when $m=5$ or $m=6$.
 By the construction of $T'$, we note that none of the paths in $T$ has order $2m-1$ or $2m-2$, and so the only possible such $T$ here has path orders $(13,13,9)$ for $m=6$ or path orders $(11,7,6)$ for $m=5$, both of which are clearly $m$-burnable.

 This completes the proof of the lemma.
\end{proof}

We are now ready to prove the main result on the burnability of path forests.

\begin{theorem}\label{n-paths}
Suppose $m\ge n\ge 3$ and let $T$ be a path forest with $n$ paths.
 If $$|T|\le m^2-(n-1)(n-2)+1-t_T\quad\mbox{and}\quad T\notin\mathcal{T}_{n,m},$$
 then $T$ is $m$-burnable.
\end{theorem}

\begin{proof}
 We shall prove the result by mathematical induction on $n$.
By Lemma~\ref{3-paths}, the base case $n=3$ follows.
 Let $n\ge 4$ and suppose the result holds for $n-1$.
 
 We will now proceed by mathematical induction on $m\ge n$.
 When $m$ equals $n$, suppose $T$ is a path forest with path orders $l_1\ge l_2\ge \cdots\ge l_n$ and $|T|\le 3n - 1 - t_T$.
 This gives $t_T = 0$ for otherwise, $|T|\ge 3(n-t_T) + 2t_T = 3n-t_T > |T|$.
 Hence, $|T|\le 3n - 1$ and $l_n = 1$.
 If $l_1\le 3$, $T$ is clearly $n$-burnable, and so we may assume that $l_1\ge 4$.
 Since $T\notin\mathcal{T}_{n,n}$ and thus $l_2\ge 2$, we also have $l_1\le 2n-1$.
 Now the path forest with $n-1$ paths obtained from $T$ by deleting its first path has order at most $3n - 1 - 4 = 3(n-1) - 2$, and so is $(n-1)$-burnable by Proposition~\ref{n-paths-3n}.
 This implies that $T$ is $n$-burnable, and so the base case $m=n$ holds.
 
 Now, consider $m\ge n+1$ and suppose the result holds for $m-1$.
 Let $T$ be a path forest with path orders $l_1\ge l_2\ge \cdots\ge l_n$ such that $|T|\le m^2 - (n-1)(n-2) +1 - t_T$ and $T\notin\mathcal{T}_{n,m}$.
 By Theorem~\ref{092021a}, we may assume that $|T| \ge m^2 - (n-1)^2 + 2$.
 Now if $l_n = 1$, we have $t_T = 0$ and $l_2\ge 2$.
 Consider the path forest $T'$ with $n-1$ paths with path orders $l_1+1, l_2+1, \ldots, l_{n-1} +1$.
 Note that 
 \begin{align*}
  |T'| &\le m^2 - (n-1)(n-2) + (n-1)\le m^2 - (n-2)^2 + 1,
 \end{align*}
 and so $T'$ is $m$-burnable by Theorem~\ref{092021a}.
 This shows that $T$ is $m$-burnable, where the first $n-1$ paths can be burned with the first $m-1$ burning sources, while the last path can be burned with the last burning source.
 So for the rest of the proof, we may assume that $l_n\ge 2$.
 
 Observe that $l_1 \ge 2m-2n+3$, as otherwise, we have $m^2 - (n-1)^2 + 2 \le |T| \le n(2m-2n + 2)$, which gives the contradiction that $(m-n)^2 + 1 \le 0$.
 We now consider a few cases based on the order of the longest path of $T$.

 \setcounter{case}{0}
 \begin{case}
  $l_1\le 2m-1$.
 \end{case}
 Consider the path forest $T'$ with path orders $l_2,l_3,\ldots,l_n$.
 Note that $t_{T'} = t_{T}$ and $|T'|\le (m-1)^2 - (n-2)(n-3) + 1 - t_{T'}$, with equality only if $l_1 = 2m-2n+3$.
 Hence, $T'$ is $(m-1)$-burnable whenever $l_1 > 2m-2n+3$ by the first induction hypothesis, and thus $T$ is $m$-burnable.
 Now for the case when $l_1 = 2m-2n+3$ and $m\ge n+2$, we can see 
 directly that $T$ $m$-burnable where each of the first $n-1$ paths can be burned with any of the first $n-1$ burning sources as $2m - 2(n-1) + 1 \ge l_1$, while the last path can be burned with the $n$th and the $(n+1)$th burning sources.
 The only case left here is when $m = n+1$ and $l_1 = 2m - 2n + 3 = 5$.
 Noting that such $T\notin\mathcal{T}_{n,n+1}$ gives $l_n \le 4$, and so $T$ is clearly $(n+1)$-burnable.

 \begin{case}
  $l_1\ge 2m$ and $l_1\neq 2m + 1$.
 \end{case}
 By considering the path forest $T'$ with path orders $l_1 - (2m-1), l_2, \ldots, l_n$, we see that $t_{T'} \le t_T$ and $|T'|\le (m-1)^2 - (n-1)(n-2) + 1 - t_{T'}$.
 Hence, $T'$ is \mbox{$(m-1)$-burnable} whenever $T'\notin\mathcal{T}_{n,m-1}$ by the second induction hypothesis, and thus $T$ is $m$-burnable.
 Since $l_n\ge 2$ and $T\notin\mathcal{T}_{n,m}$, we see that for $T'\in\mathcal{T}_{n,m-1}$, we must have % $|T| = m^2 - (n-1)(n-2)+1 -t_T$ and 
 $l_1 = 2m+4$ or $l_1 = 2m+2$.
 Consider such a path forest $T$.

 If $l_1 = 2m+4$, we have $l_2 = (m-1)^2 - (n-1)(n+3) + 1$ and $l_3 = l_4 = \cdots = l_n = 5$.
 Since $T\notin\mathcal{T}_{n,m}$, we also have $l_2\ge 6$.
 Now observe that $2m+4 = l_1\ge l_2 = (m-n-2)(m+n)+5\ge 6$ is only possible when $m = n+3$.
 So we have $l_1 = 2m+4 = (2m-3) + 7$, $l_2=2m+2 = (2m-1) + 3$, and it is straightforward to verify that $T$ is $(n+3)$-burnable using only the first $n+2$ burning sources.

 If $l_1 = 2m+2$, we have $l_2 = (m-1)^2 - n^2 + 2$ and $3\ge l_3, l_4, \ldots, l_n \ge 2$.
 Since $T\notin\mathcal{T}_{n,m}$, we also have $l_2\ge 4$.
 Now observe that $2m+2 = l_1\ge l_2 = (m-n-1)(m+n-1)+2\ge 4$ is only possible when $m = n+2$.
 So we have $l_1 = 2m+2 = (2m-3) + 5$, $l_2=2m-1$, and it is straightforward to verify that $T$ is $(n+2)$-burnable using only the first $n+1$ burning sources.

 \begin{case}
  $l_1 = 2m+1$.
 \end{case}
 
 We first observe that in this case, we have $m\ge n+2$, since if $m = n+1$, then
 $$|T|\ge (2n+3) + 3(n-t_T-1) + 2t_T = m^2 - (n-1)(n-2) + 1 - t_T\ge |T|,$$
 which would only be possible if $T\in\mathcal{T}_{n,n+1}$.
 We now proceed with the following claim using similar arguments as in  the previous cases.
 
 \begin{claim}
  If $2m-2n+3\le l_{i_0}\le 2m$ for some $2\le i_0\le n$, then $T$ is $m$-burnable.
 \end{claim}

(Proof of Claim)  If there is such an $i_0$, we consider the path forest $T'$ obtained from $T$ by deleting $\min\{2m-1, l_{i_0}\}$ vertices from its $i_0$ path.
  Then either $T'$ is a path forest consisting of $n$ paths where $ \vert T'\vert\leq  (m-1)^2 - (n-1)(n-2) + 1$ and the shortest path of $T'$ has one vertex, or $T'$ is a path forest consisting of $n-1$ paths and $\vert T'\vert \leq   (m-1)^2 - (n-2)(n-3) + 1 - t_{T'}$ with $t_{T'}= t_T$.
  
  For the former, noting that $l_n\ge 2$ and so $T'\notin\mathcal{T}_{n,m-1}$, we see that $T'$ is $(m-1)$-burnable by the second induction hypothesis.
  For the latter, we further delete $2m-3$ vertices from the longest path of $T'$ (it is also the first path of $T$) to obtain a path forest $T''$ consisting of $n-1$ paths where one of its path has order $l_1 - (2m-3)=4$ and $t_{T''} = t_T$. 
  It can be verified now that $T''\notin\mathcal{T}_{n-1,m-2}$ and so noting that $|T''|\le (m-2)^2 - (n-2)(n-3) + 1 - t_{T''}$, we see that $T''$ is $(m-2)$-burnable by the first induction hypothesis.
  This in turn implies that $T'$ is $(m-1)$-burnable.
  
  In both cases, $T'$ is $(m-1)$-burnable. It follows that $T$ is $m$-burnable, completing the proof of the claim.
 
 With the claim above, it remains to consider $T$ with path orders
 $$ \underbrace{2m+1,2m+1,\dotsc, 2m+1}_{k \text{ times for some } k\geq 1 },  \underbrace{l_{k+1}, l_{k+2},\dotsc, l_n}_{\text{each } \leq 2m-2n+2}.$$
 We now show that $m \ge n+k+1$.
 Recalling that $m\ge n+2$, this is true for when $k=1$.
 For $k\ge 2$, we have
 \begin{align*}
  m^2 - (n-1)(n-2) + 1 - t_T\ge |T| &\ge k(2m+1) + 3(n-k) - t_T,
 \end{align*}
 which gives $$(m-k)^2\ge n^2 + (k-1)^2 > n^2,$$
 and thus $m>n+k$.
 
 Finally, we show directly that $T$ is $m$-burnable as follows. 
 Noting that $2m + 1 \le (2m - 2i + 1) + (2i + 1)$, each of the first $k$ paths of order $2m+1$ is burned with the $i$th and the $(m-i)$th burning sources for any $1\le i\le k$.
 This way, the $(k+1)$th until the $n$th burning sources and the last burning source are yet to be used.
 As the $(n-1)$th burning source can burn $2m - 2n + 3$ vertices, while the $n$th and the last burning sources can collectively burn $2m-2n+2$ vertices, we see that the last $n-k$ paths of $T$, each of which has order at most $2m-2n+2$, can be burned with these $n-k+1$ burning sources.
 Therefore, $T$ is $m$-burnable.
 
 This completes the proof of the theorem.
\end{proof}

% This is the edited sections on double spiders
%
%
%
%
%
%
\section{Double Spiders}\label{double}

We study the burnability of double spiders in this section.
We start with \mbox{$n$-double} spiders of order $m^2+n-2$ where $m>n$, showing that such double spiders are $m$-burnable.
This verifies Conjecture~\ref{Tree Conjecture} for double spiders.

Clearly, a $2$-double spider is either a path or a $3$-spider, and so it is $m$-burnable when its order is at most $m^2$.
The following easy lemma asserts that this can in fact be done in such a way that the heads are burned in the earlier rounds.
This lemma will serve as the base case of our first result on double spiders.

\begin{lemma}\label{double_n=2}
 Let $m\geq 3$.
 Then every $2$-double spider of order at most $m^2$ is \mbox{$m$-burnable}.
 Furthermore, if the shortest arm has length $l$, then 
 there is a way to burn the double spider in $m$ rounds such that at least $\min\{l,m-2\}$ rounds still remain after both its heads are burned.
 \end{lemma}

\begin{proof}
 Suppose $T$ is a $2$-double spider with arms of lengths $l_1\geq l_2$ and its order is at most $m^2$.
 For convenience, we may as well assume that $|T| = m^2$.
 Regardless of whether $T$ is a path or a $3$-spider, the following burning strategy works.

 If $l_2\le m-2$, then we put the first burning source at the $(m-2-l_2)$th vertex on the first arm. In $m$ rounds, this first burning source would burn the entire second arm, and clearly  after both the heads are burned by this first burning source, there are still $l_2$ rounds left.
 The remaining vertices unburned by the first burning source form a path of order $(m-1)^2$, and so is clearly $(m-1)$-burnable.

 Suppose now that $l_2\ge m-1$.
 Then we put the first burning source at one of the heads so that the remaining at most $(m-1)^2$ vertices unburned by the first burning source form one of the following graphs $T'$.
 Either $T'$ is a path (possible when $l_2 = m-1$), or $T'$ is a path forest with two paths such that its path orders are not $(m-1)^2 - 2$ and $2$.
 Therefore, $T'$ is $(m-1)$-burnable, and thus the lemma follows.
\end{proof}

In the above lemma, we see that when both arms of a $2$-double spider are just long enough (at least $m-2$), we have an optimal burning sequence starting at one of its heads.
But just like burning spiders, this is not always the case for double spiders with more than two arms.
Consider the $3$-spider with arms of lengths $5$, $5$, and $6$.
Then any optimal burning sequence (of length four) has to start from the first vertex next to the head on one of the arms of length five.
Hence, if we regard this as a $3$-double spider with three arms of equal lengths, we see that for any optimal burning sequence, after both heads are burned, there is only one round left.
So we may not always start the burning sequence at one of the heads for an $n$-double spider with $n\geq 3$ even when all of its arms are just long enough, slightly different from that of Lemma~\ref{double_n=2}.
Just like the corresponding result on spiders in \cite{tan2020graph}, extending Theorem~\ref{double_spiders} to keep track on when the heads are burned in a double spider simplifies the proof.

\begin{theorem}\label{double_m>n}
 Suppose $m>n\geq 2$. Then every $n$-double spider of order at most $m^2+n-2$ is  $m$-burnable. Furthermore, if  the shortest arm has length $l$, then there is a way to burn the double spider in $m$ rounds such that at least $\min\{l,m-3\}$ rounds still remain after both its heads are burned.
\end{theorem}

\begin{proof}
 We prove by mathematical induction on $n$.
 By Lemma~\ref{double_n=2}, the base case $n=2$ follows.
 Suppose $n\ge 3$ and the result holds for $n-1$.
 Let $m>n$ and suppose $T$ is an $n$-double spider with arms of lengths $l_1\ge l_2\ge \cdots\ge l_n$   and $\vert T\vert \leq m^2+n-2$.
 We may as well assume that $|T|=m^2+n-2$.
 We start with the following claim that deals with the case when the shortest arm is not too long.
 \begin{claim}
  We may assume that $l_n\ge m$.
 \end{claim}
 
 (Proof of Claim)
  First, we suppose $l_n\le m-3$.
  Let $T'$ be the $(n-1)$-double spider obtained from $T$ by deleting its $n$th arm.
  Since $|T'| \le m^2+(n-1)-2$, we see from induction hypothesis that there is a way to burn  $T'$ in $m$ rounds such that at least $\min\{l_{n-1}, m-3\}\ge l_n$ rounds remain  after both its heads are burned.
  Clearly, $T$ can be burned by the same burning sequence.

  Now, we suppose $l_n$ is either $m-2$ or $m-1$.
  If $l_n=m-1$, we put the first source at the head where the $n$th arm is joined to, while if $l_n= m-2$, we put the first source at the head where some arm of length at least $m-1$ is joined to.
  In either case, in $m$ rounds, the first source would burn at least $n(m-2)+1$ vertices from the $n$ arms as well as both the heads, and the rest of the vertices form a path forest $T'$ with at most $n-1$ paths.
  Note that
  \begin{align*}
    |T'| &\leq  m^2+n-2-n(m-2)-3\\
		 &=(m-1)^2+3n-m(n-2)-6\\
		 &\le (m-1)^2+3n-(n^2-n-2)-6  &\text{(because $m\geq n+1$)}   \\
		 &= (m-1)^2-(n-2)^2.
  \end{align*}
  By Theorem~\ref{092021a}, $T'$ is $(m-1)$-burnable and therefore $T$ can be burned in $m$ rounds starting at one of its heads. This completes the proof of the claim.

 With the above claim, we now have $l_i\ge m$ for all $1\le i\le n$.
 We will show that $T$ is $m$-burnable, and that in most cases, $T$ has an optimal burning sequence starting at one of its heads, while in the remaining few cases, the optimal burning sequence starts at the first vertex of an arm next to the head it is joined to.
 With such burning sequences, at least $m-3$ rounds are left after both heads of $T$ are burned, which would complete the induction step.

 Our burning strategy for every case is similar.
 We first put the burning source at one of the heads or adjacent to them (i.e.~the first vertex of an arm), and we observe that the remaining vertices unburned by this first burning source (after $m$ rounds) form a path forest with at most $n$ paths, which we will denote as $T'$.
 Finally, by showing that $T'$ is $(m-1)$-burnable using Theorem~\ref{092021a} and Proposition~\ref{n-paths-linear}, we see that $T$ is $m$-burnable.

 We now consider these cases.

 \setcounter{case}{0}

 \begin{case}
  $m\geq n+3$.
 \end{case}
 The first source is placed at the head of $T$ where at least two arms are joined to it (possible as $n\geq 3$).
 This first source would burn at least $m-1$ vertices from the two aforementioned arms in $m$ rounds, and at least $m-2$ vertices from each of the other arms.
 As before, we see that
 \begin{align*}
  |T'| &\leq m^2+n-4-2(m-1)-(n-2)(m-2)\\
       &= (m-1)^2+n-(n-2)(m-2)-3\\
       &\leq (m-1)^2+n-(n^2-n-2)-3 &\text{(because $m\geq n+3$)}   \\
       &= (m-1)^2-(n-1)^2.
 \end{align*}
Hence,  $T'$ is $(m-1)$-burnable by Theorem~\ref{092021a}.

 \begin{case}
  $m=n+1$.
 \end{case}
 If $l_n = m$, then we put the first source this time at the head where the $n$th arm is joined to.
 Letting the path orders of $T'$ be $l_1',l_2',\ldots,l_n'$, we see that $l_i'\le l_i - (m-2) =: l_i''$ for $1\le i\le n-1$ and $l_n' = 1 =: l_n''$.
 Now, since $l_i''\ge 2$ for $1\le i\le n-1$, $l_n''=1$, and
 $$\sum_{i=1}^n l_i'' = \left(\sum_{i=1}^n l_i\right) - n(m-2) - 1 = m^2 + n - 4 - n(m-2) - 1 = 4n - 4,$$
 Proposition~\ref{n-paths-linear}\eqref{n-paths-4n} implies that the path forest with path orders $l_1'',l_2'',\ldots,l_n''$ is $n$-burnable.
 This in turn implies that $T'$ is also $n$-burnable.

 If $l_n\ge m+1$, then it is straightforward that $T$ has at least three arms of length exactly $m+1$.
 Indeed, letting $k$ to be the number of arms of $T$ of length $m+1$, we see that
 $$ n^2 + 3n - 3 = m^2 + n - 4 \ge (n-k)(m+2) + k(m+1) = n(n+3)-k,$$
 which clearly implies $k\ge 3$.
 So we may suppose the $(n-1)$th and the $n$th arms are joined to the same head, and we put the first source at the vertex on the $n$th arm next to the head.
 Observe that this first burning source would burn exactly $m-2$ vertices and $m$ vertices of the last two arms respectively, while at least $m-3$ vertices from each of the other arms.
 As above, to show that $T'$ is $n$-burnable, it suffices to show that the path forest with path orders$$l_1-(m-3),l_2-(m-3),\ldots,l_{n-2}-(m-3), 3, 1$$ is $n$-burnable.
 Since this path forest has order $m^2 + n - 4 - n(m-3) - 4 = 5n-7$, it is $n$-burnable by Proposition~\ref{n-paths-linear}\eqref{n-paths-5na}.

 \begin{case}
  $m=n+2$.
 \end{case}
 If $l_n=m$, the first burning source is put again at the vertex on the $n$th arm next to the head.
 Observe that this first burning source would burn the $n$th arm completely and at least $m-3$ vertices from each of the other arms.
 So $T'$ has $n-1$ paths and
 $$ |T'| \le m^2 + n - 4 - (n-1)(m-3) - m = 6n - 3= (n+1)^2 - (n-2)^2,$$
 which is $(n+1)$-burnable by Theorem~\ref{092021a}.

 For the final case when $l_n\ge m+1$, we first note that $l_1\ge n+5 = m+3$ as $|T| = m^2+n - 2 = n^2+5n+2$.
 We choose the head of $T$ to which the first arm is joined as where to put our first burning source.
 Observe that this first source would burn the first $m-1$ vertices from the first arm, and at least $m-2$ vertices from each of the other arms.
 As in Case 2, we want to show that the path forest with path orders
 $$l_1 - (m-1), l_2 - (m-2), l_3-(m-2), \ldots, l_n - (m-2)$$
 is $(n+1)$-burnable.
 Noticing that this path forest has order $n^2+5n - n(m-2) - 1 = 5n -1$ and each of its paths has order at least three, we see that it is $(n+1)$-burnable by Proposition~\ref{n-paths-linear}\eqref{n-paths-5nb}.
 This of course implies that $T'$ is $(n+1)$-burnable.

 From these cases, the proof of the theorem is now complete.
\end{proof}

For the remaining of the section, we consider $n$-double spiders of order $m^2+n-2$ where $m\le n$, proving the corresponding result on spiders in \cite{tan2020graph}.
In this case, not all such double spiders are $m$-burnable, as observed by the following simple lemma.

\begin{lemma}
 Let $m\geq 2$. If $T$ is an $m$-double spider such that the distance between any two leaves is at least $2m$, then $T$ is not $m$-burnable.
\end{lemma}

\begin{proof}
 Suppose $T$ is an $m$-double spider as in the lemma.
 The unique path joining any two leaves has order at least $2m+1$. Consider a burning process that takes $m$ rounds. Suppose a source is put on this path and it burns one of the two leaves when the burning  process ends.
 Then this burning source would not burn at least two vertices at the opposite end of the path.
  It follows that in a burning process of $m$ rounds where $m-1$ leaves are burned by the first $m-1$ sources at the end of the process, the last remaining leaf together with its neighbor would not be burned by these $m-1$ burning sources.
 However, the last burning source burns only one vertex, and thus $T$ would not be completely burned at the end of the process.
 Therefore, $T$ is not $m$-burnable.
\end{proof}

The following theorem completes our results on the burnability of double spiders.

\begin{theorem}\label{double_m<n}
 Suppose $n\ge m\ge 3$ and let $T$  be an $n$-double spider such that $\vert T\vert \leq m^2+n-2$. Then $T$ is $m$-burnable unless it includes, as a subgraph, an $m$-double spider such that the distance between any two leaves is at least $2m$.
\end{theorem}

\begin{proof}
 Let $T$ be an $n$-double spider of order $m^2+n-2$ and suppose it does not contain an $m$-double spider such that the distance between any two leaves is at least $2m$.
 Letting $l_1\ge l_2\ge \cdots \ge l_n$ be the lengths of the arms of $T$ as usual, we note that it must be the case that $l_m\le m-1$.
 We shall assume that $l_1\ge m-1$, for otherwise $T$ is clearly $m$-burnable.
 Now, we consider a few cases based on $l_m$.
 
 \setcounter{case}{0}
 
 \begin{case}
  $l_m\le m-3$.
 \end{case}
  Note that $m\ge 4$.
 Let $T'$ be the $(m-1)$-double spider obtained from $T$ by deleting its last $n-m+1$ arms. Hence, 
 $$|T'| \le m^2 + n - 2 -  (n-m+1) = m^2 + (m-1) - 2.$$
 So by Theorem~\ref{double_m>n}, there is a way to burn $T'$  in $m$ rounds such  that after both of its heads are burned, at least $\min\{l_{m-1}, m-3\}$ rounds remain.
 Since $\min\{ l_{m-1}, m-3\}  \ge l_m$, it follows that the same burning sequence will burn $T$ as well.
 Thus, $T$ is $m$-burnable.
 
 \begin{case}
  $l_m = m-2$.
 \end{case}
 Note that in this case, we must have $l_{m-1}\le m$, for otherwise $|T| \ge (m-1)(m+1) + (m-2) + n - m + 2 = m^2 + n - 1 > |T|$.
 To show that $T$ is $m$-burnable, we will place the first burning source at one of its heads.
 Observe that in $m$ rounds, this burning source would clearly burn the last $n-m+1$ arms of $T$ completely and at least $m-2$ vertices from each of the first $m-1$ arms.
 By appropriately choosing the head for the first burning source as follows, it can be guaranteed that one extra vertex would be burned by this burning source.
 If $l_{m-1} = m-2$, we choose the head where the first arm is joined to, and if $l_{m-1}\ge m-1$, we choose the head where the $(m-1)$th arm is joined to.
 This way, we also have that at most one vertex of the $(m-1)$th arm is unburned by the first burning source.

 Letting $T'$ be the path forest formed by the remaining vertices unburned by the first burning source, we see that $T'$ has at most $m-1$ paths, and if it has exactly $m-1$ paths (only possible if $l_{m-1} = m$), its shortest path has order one.
 Of course, $T'$ contains, as a spanning subgraph, a path forest with exactly $m-1$ paths where its shortest path has order one.
 Since $$|T'|\le m^2 + n - 4  - m(m-2) - 1 - (n-m) = 3m - 5 = 3(m-1) - 2,$$
 this spanning subgraph, and thus $T'$, is $(m-1)$-burnable by Proposition~\ref{n-paths-3n}.
 It follows that $T$ is $m$-burnable.

 \begin{case}
  $l_m = m-1$.
 \end{case}
  We will show in this case that $T$ is $m$-burnable by putting the first source at the head where the $m$th arm is joined to.
  As in Case 2, the remaining vertices unburned by this first burning source (after $m$ rounds) form a path forest $T'$ of order at most $3m - 5$.
  We claim that $T'$ has at most $m-2$ paths, or $T'$ has exactly $m-1$ paths with the shortest path having order one, and so again like in Case 2, $T'$ is $(m-1)$-burnable.

  If $l_{m+1}\le m-2$ (or if $m = n$), we note that among the first $m-1$ arms of $T$, either some arm that joins to the head where the $m$th arm is joined to has length at most $m$, or some arm that joins to the other head has length at most $m-1$.
  This is due to the property of $T$ that it does not contain an $m$-double spider such that any two leaves have distance at least $2m$.
  It is now clear that $T'$ is as claimed.

  If $l_{m+1} = m - 1$ and $l_{m+2}\le m-2$ (or if $l_{m+1} = m-1$ and $n=m+1$), we note by simple counting that we must also have $l_{m-1} = m-1$.
  So we may assume both the $m$th and the $(m+1)$th arms of $T$ are joined to the same head, and again, we see that $T'$ is as claimed.

  Finally if $l_{m+2} = m-1$, it must be the case that the first $m+2$ arms of $T$ each has length exactly $m - 1$, while the last $n - m - 2$ arms each has length exactly one.
  So we may assume the $m$th arm of $T$ is joined to the head where the majority of the longer arms are joined to.
  Since $m+2\ge 5$, $T'$ is as claimed.
  In fact, $T'$ consists of at most $\lfloor\frac{m}{2}\rfloor+1$ independent vertices.

  From the above cases, the proof of the theorem is complete.
\end{proof}

\section{Path Forests with Sufficiently Long Paths}\label{path-forests-long}

This section  aims to show that a path forest of order $m^2$ with a sufficiently long shortest path is always $m$-burnable, as stated in Theorem~\ref{long-path-forests}.
For this purpose, we may assume in this section that every path forest under our consideration has order $m^2$ for some $m\in\mathbb{N}$.
Additionally, we often identify a path forest with $n$ paths with the $n$-tuple listing of its path orders $(l_1,l_2, \ldots, l_n)$ and we \emph{do not} always assume $l_1\ge l_2\ge \cdots \ge l_n$.
Generalizing from the context of burning connected graphs, we can also say that a path forest of order $m^2$ is \emph{well-burnable} if it is $m$-burnable.
When a path forest is not well-burnable, we say it is \emph{deficient}.
We start by defining a relation for comparing `successive' deficient path forests.

\begin{definition}
 Let $n\geq 2$.
 Suppose the deficient path forests $T =(l_1,l_2, \ldots, l_n)$ and $T'=(l'_1,l'_2, \ldots, l'_n)$ have orders $m^2$ and $(m+1)^2$, respectively.
 We denote by $T\prec T'$ if there exists $1\leq i\leq n$ such that $l'_i= l_i+(2m+1)$ and $l'_j=l_j$ for all $j\neq i$.
 We say that $T'$ is obtained from $T$ by \emph{$\prec$-extension at the $i$-th component}, and we also write $T'\succ T$.
\end{definition}

Starting from a deficient path forest $T$, there are two possibilities for successive $\prec$-extensions, either
\begin{enumerate}
 \item there exists an infinite chain $T=T_0 \prec T_1 \prec T_2 \prec \dotsb$ starting from $T$; or
 \item no such infinite $\prec$-chains exist starting from $T$.
\end{enumerate}  
For example, $(7,2,2)\prec (12,2,2) \prec (21,2,2) \prec (34,2,2)\prec \dotsb $ is an infinite chain starting from the deficient path forest $(7,2,2)$, whereas $(17,15,4)\prec (17,15,17)\prec (17,30,17)$, $(17,15,4)\prec (17,28,4)\prec (17,43,4)$ and $(17,15,4)\prec (30,15,4)\prec (30,30,4)$ are chains starting from the deficient path forest $(17,15,4)$ which cannot be extended further.

The following notation is convenient for stating and proving our results here.

\begin{notation}
For $n,L\in\mathbb{N}$, let $\mathcal{H}(n,L)$ be the set of all path forests with $n$ paths where each of its paths has order at least $L$.
The set of all deficient path forests in $\mathcal{H}(n,l)$ is denoted by $\mathcal{H}_{\text{def}}(n,l)$.
\end{notation}

\begin{theorem}\label{110322c}
 For each $n\geq 2$, the following are equivalent.
 \begin{enumerate}
  \item\label{first} There exists $L\in\mathbb{N}$ such that every $T\in \mathcal{H}(n,L) $ is well-burnable.
  \item\label{second} There exists $L'\in\mathbb{N}$ such that there are no infinite $\prec$-chains starting from a path forest $T$ in $\mathcal{H}_{\text{def}}(n,L')$.
 \end{enumerate}
\end{theorem}

\begin{proof}
Fix the number of paths $n$.
Assume \eqref{first} holds and let $L$ has the property that every $T\in\mathcal{H}(n,L) $ is well-burnable.
Then \eqref{second} vacuously follows as $\mathcal{H}_{\text{def}}(n,L)$ is empty and so we can pick $L'=L$.
Conversely, assume \eqref{second} holds and let $L'$ has the property that there are no infinite $\prec$-chains starting from a path forest in $\mathcal{H}_{\text{def}}(n,L')$.
For each $T \in \mathcal{H}_{\text{def}}(n,L') $, we see from the Konig's Tree Lemma that the tree rooted at $T$ induced by the relation $\prec$ is finite.
Indeed, every vertex of such a rooted tree has at most $n$ descendants, and there is no infinite path by our assumption.

Now, observe that for sufficiently large $m$, the longest path of a path forest of order $m^2$ has order at least $L' + (2m-1)$.
Hence, if $T \in \mathcal{H}_{\text{def}}(n,L')$ has large enough order, then there exists $T' \in \mathcal{H}_{\text{def}}(n,L')$ such that $T'\prec T$.
This implies that there is some $M\in \mathbb{N}$ such that for every $T \in \mathcal{H}_{\text{def}}(n,L')$ of order $m^2$ where $m\ge M$, there is a decreasing chain $T =T_m \succ T_{m-1}\succ \cdots \succ T_M$ with $T_i \in \mathcal{H}_{\text{def}}(n,L')$ and $|T_i| = i^2$  for every $M\leq i \leq m$; in other words, $T$ is a vertex belonging to the finite tree rooted at $T_M$.
Since there are finitely many $T\in \mathcal{H}_{\text{def}}(n,L')$ of order at most $M^2$, the set $\mathcal{H}_{\text{def}}(n,L')$ is finite.
Pick a large enough integer $L$ such that $\mathcal{H}(n, L) \cap \mathcal{H}_{\text{def}}(n,L') =\emptyset  $.
It follows that every $T\in \mathcal{H}(n, L)$ is well-burnable, completing the proof of the theorem.
\end{proof}

The above theorem provides a reduction of our main result to an equivalent result that we will actually prove, namely \eqref{second} in Theorem~\ref{110322c}.
To show that there are no such infinite $\prec$-chains, we will need a series of technical lemmas, starting with Lemma~\ref{090322a} and Lemma~\ref{080322a}, which essentially imply that a deficient path forest with sufficiently long paths cannot be $\prec$-extended infinitely often only at the same one component.

\begin{lemma}\label{090322a}
 Let $n\in \mathbb{N}$, and suppose $l_1,l_2,\ldots, l_n$ are even numbers, each at least $8n$.
 Then there is some $m$ such that the path forest $(l_1, l_2, \ldots, l_n, m^2-\sum_{i=1}^n l_i)$ has an optimal burning sequence where the last $n$ burning sources are not used to burn the first $n$ paths.
\end{lemma}

\begin{proof}
 We may assume that $l_1\ge l_2\ge \cdots \ge l_n$.
 For each $1\le i \le n$, let $t_i$ be the integer such that $l_i =  (2t_i-1) + (2(n+i)-1)$.
 We can take $m=t_1$.
 It suffices to verify that $m = t_1>t_2>\cdots > t_n>2n$, as this would imply that the path forest $(l_1, l_2, \ldots, l_n, m^2-\sum_{i=1}^n l_i)$ has an optimal burning sequence where the first $n$ paths are burned using the respective $2n$ distinct burning sources that exclude the last $n$ burning sources.
 Clearly, $t_i>t_{i+1}$ as $l_{i}\ge l_{i+1}$ for every $1\le i\le n-1$, and the inequality $t_n>2n$ follows from the condition $l_n\ge 8n$.
\end{proof}

\begin{lemma}\label{080322a}
 Let $n\in \mathbb{N}$.
 There exists $L\in \mathbb{N}$ such that whenever $l_1,l_2,\ldots, l_n$ are integers at least $L$, the path forest $(l_1, l_2, \ldots, l_n, m^2-\sum_{i=1}^n l_i)$ is $m$-burnable for some $m\in\mathbb{N}$.
\end{lemma}

\begin{proof}
 For a fixed $n$, we can take $L=10n-1$.
 Let $l_1,l_2,\ldots,l_n$ be some given integers, each at least $10n-1$.
 Without loss of generality, suppose $l_1, l_2, \ldots, l_k$ are odd while $l_{k+1}, l_{k+2}, \dotsc, l_n$ are even, where $0\leq k \leq n$.
 Let $l'_i= l_i-(2i-1)$ for $0< i \le k$ and $l'_i= l_i$ for $k<i\leq n$.
 Note that $l'_i \geq 8n$ for every $1\leq i\leq n$.
 By Lemma~\ref{090322a}, for some $m\in \mathbb{N}$, the path forest $(l'_1, l'_2, \dotsc, l'_n, m^2-\sum_{i=1}^n l'_i)$ has an optimal burning sequence where the last $n$ burning sources are not used to burn the first $n$ paths.
 Therefore, by reallocating the last $k$ burning sources accordingly, it follows that the path forest $(l_1, l_2, \dotsc, l_n, m^2-\sum_{i=1}^n l_i)$ is $m$-burnable.
\end{proof}

We remark that with some careful analysis, one can show that for $n=2$, the least $L$ that satisfies the property stated in Lemma~\ref{080322a} is $L=8$.
However, as we will see later, Lemma~\ref{080322a} will not give any quantifiable bound for our main result, and so we do not attempt to optimize the bound of $L$ in the lemma.

The main technical lemma needed for our main result is Lemma~\ref{070322b}, which allows us to show that if we go very far along a $\prec$-chain starting from a path forest with $n$ paths, we will have the flexibility to allocate burning sources accordingly so that any of the $n$ next potential extensions is well-burnable, and hence terminating the $\prec$-chain.
The following lemma is essentially the base case for the induction proof of Lemma~\ref{070322b}.
%This lemma can be interpreted informally as follows: when the path lengths are sufficiently long depending on the number of paths, then given all the firepowers $1,3,5,7,\dotsc$ at our disposal, we can burn the path forest using some of the firepowers without any wastage.

\begin{lemma}\label{070322a}
 Suppose $\langle a_i\rangle_{i=1}^\infty$ and $\langle b_i\rangle_{i=1}^\infty$ are strictly increasing sequences of odd integers without common terms such that $a_i+2\in \{ b_j \mid j\in \mathbb{N}\}$ for infinitely many $i$.
 Then for any given integer $x$, there exist indices $N_1$ and $N_2$ such that the set $\{a_i \mid 1\le i\le N_1\}\cup  \{b_i \mid 1\le i\le N_2\}$ can be partitioned into two sets $C\cup D$ with the property that $$\summation(C)=  \sum_{i=1}^{N_1} a_i +x \quad\mbox{ and }\quad \summation(D)= \sum_{i=1}^{N_2} b_i-x.$$
\end{lemma}

\begin{proof}
 The result is trivial for $x=0$.
 So we suppose first that $x=2k$ for some positive integer $k$.
 Let $i_1, i_2, \ldots, i_k$ enumerate the least $k$ terms of $\langle a_i\rangle_{i=1}^\infty$ with $a_{i}+2 \in \{ b_j \mid j\in \mathbb{N}\}$.
 Let $N_1= i_k$ and let $N_2$ be the index in $\langle b_i\rangle_{i=1}^\infty$ such that $b_{N_2}= a_{i_k}+2$.
 By taking the partition $\{a_i \mid 1\le i\le N_1\}\cup  \{b_i \mid 1\le i\le N_2\} = C\cup D$, where
 \begin{align*}
  C &= \left(\{ a_i \mid 1\le i\le N_1 \}\backslash \{ a_{i_j}\mid  1\leq j\leq k \} \right) \cup   \{ a_{i_j}+2\mid  1\le j\le k \} \text{ and }  \\
  D &=\left(\{ b_i \mid 1\le i\le N_2 \}\backslash \{ a_{i_j}+2\mid  1\leq j\leq k \} \right)\cup   \{ a_{i_j}\mid  1\le j\le k \},
 \end{align*}
 it is easy to see that $C$ and $D$ have the required property.

 Now, suppose $x$ is odd or negative.
 Choose $l$ such that $x + \sum_{i=1}^l a_i$ is even and nonnegative. (Clearly, such an $l$ exists.)
 Applying what we have proved above on the sequences  $\langle a_i\rangle_{i=l+1}^\infty$ and $\langle b_i\rangle_{i=1}^\infty$, it follows that for some $N_1$ and $N_2$, the set $\{ a_i \mid l+1\leq i\leq N_1 \}\cup  \{ b_i \mid 1\leq i\leq N_2\}$ can be partitioned into $C'\cup D'$ with
 $\summation(C')=  \sum_{i=l+1}^{N_1} a_i + (x+ \sum_{i=1}^l a_i) $ and $\summation(D')= \sum_{i=1}^{N_2} b_i- (x+\sum_{i=1}^l a_i)$.
 It is now straightforward to see that the sets $C := C'$ and $D := D' \cup \{a_i \mid 1\le i\le l\}$ have the required property.
%Take $C'= C$ and $D'= D	\cup  \{ a_i \mid 1\leq i\leq l\}$.
%Then $C'$ and $D'$ have the required property.
\end{proof}

\begin{lemma}\label{070322b}
 Let $n\geq 2$.
 Suppose $\langle a_{1,i}\rangle_{i=1}^\infty, \langle a_{2,i}\rangle_{i=1}^\infty, \ldots, \langle a_{n,i}\rangle_{i=1}^\infty$ are $n$ strictly increasing sequences of odd integers without common terms such that their union $\bigcup_{j=1}^n\{a_{j,i} \mid i\in \mathbb{N}\}$ is the set of all odd integers at least $\min\{a_{1,1}, a_{2,1}, \ldots, a_{n,1}\}$.
 Then for any $n$ integers $x_1, x_2,\ldots, x_{n}$ such that $\sum_{j=1}^n x_j=0$, there exist indices $N_1, N_2, \ldots, N_n$ such that the set $\bigcup_{j=1}^n \{ a_{j,i} \mid 1\le i\le N_j \}$ can be partitioned into $C_1\cup C_2\cup \cdots \cup C_n$ with the property that for every $1\le j\le n$,
 $$\summation(C_j) =  \sum_{i=1}^{N_j} a_{j,i} +x_j.$$
\end{lemma}

\begin{proof}
 We prove by mathematical induction on $n$.
 The case $n=2$ follows from Lemma~\ref{070322a}, as it can be deduced in this case that $a_{1,i}+2\in \{ a_{2,i'} \mid i'\in \mathbb{N}\}$ for infinitely many $i$.
 For the induction step, suppose $\langle a_{1,i}\rangle_{i=1}^\infty, \langle a_{2,i}\rangle_{i=1}^\infty, \ldots, \langle a_{n+1,i}\rangle_{i=1}^\infty$ are $n+1$ strictly increasing sequences of odd integers without common terms such that the union of their terms $\bigcup_{j=1}^{n+1}\{ a_{j,i} \mid i\in \mathbb{N} \}$ is the set of all odd integers at least $\min\{a_{1,1}, a_{2,1}, \ldots, a_{n+1,1}\}$.
 We first note that there are infinitely many elements $z$ in $\bigcup_{j=1}^n \{ a_{j,i} \mid i\in \mathbb{N} \}$ such that $z+2\in \{ a_{n+1,i} \mid i\in \mathbb{N} \} $, for otherwise $\{ a_{n+1,i} \mid i\in \mathbb{N} \}$ is either finite or cofinite.
 By the infinite pigeonhole principle, we may suppose that $a_{n,i}+2\in \{ a_{n+1,i'} \mid i'\in \mathbb{N}\}$ for infinitely many $i$.

 Let integers $x_1, x_2,\ldots, x_{n+1}$ with $\sum_{j=1}^{n+1} x_j=0$ be given.
 We merge $\langle a_{n,i}\rangle_{i=1}^\infty $  and $\langle a_{n+1,i}\rangle_{i=1}^\infty$ into one strictly increasing sequence $\langle b_{i}\rangle_{i=1}^\infty$.
 Consider the $n$ sequences $\langle a_{1,i}\rangle_{i=1}^\infty  , \langle a_{2,i}\rangle_{i=1}^\infty, \ldots, \langle a_{n-1,i}\rangle_{i=1}^\infty, \langle b_i\rangle_{i=1}^\infty$ and the $n$ integers $x_1, x_2,\ldots, x_{n-1}, x_n+ x_{n+1}$.
 By induction hypothesis, for some indices $N_1, N_2, \ldots, N_n$, there exists a partition of
 $$\bigcup_{j=1}^{n-1} \{ a_{j,i} \mid 1\le i\le N_j \} \cup \{b_i\mid 1\le i\le N_n\} = C_1\cup C_2\cup \cdots \cup C_n$$ with the property that for each  $1\le j\le n-1$,
 $$\summation(C_j)= \sum_{i=1}^{N_j} a_{j,i} + x_j, \quad\text{and}\quad
 \summation(C_n)= \sum_{i=1}^{N_n} b_{i} + x_n+x_{n+1}.$$
 By construction, we have $\{b_i\mid 1\le i\le N_n\}= \{ a_{n,i}\mid 1\le i\le M_1\}\cup \{ a_{n+1,i}\mid 1\le i\le M_2\}$ for some $M_1$ and $M_2$.
 
 Now, consider the two increasing sequences $\langle a_{n,i}\rangle_{i=M_1+1}^\infty$ and $\langle a_{n+1,i}\rangle_{i=M_2+1}^\infty$ and set $x= -\sum_{i=1}^{M_2} a_{n+1,i}  -x_{n+1}$.
 By Lemma~\ref{070322a}, for some indices $N'_n$ and $N'_{n+1}$, 
 there exists a partition of
 $$\{a_{n,i} \mid M_1+1\le i\le N'_n\}\cup \{a_{n+1,i} \mid M_2+1\le i\le N'_{n+1}\} = D_1\cup D_2$$
 with the property that 
 $$\summation(D_1)=  \sum_{i=M_1+1}^{N'_n} a_{n,i} +x \quad \text{and}\quad \summation(D_2)= \sum_{i=M_2+1}^{N'_{n+1}} a_{n+1,i}-x.$$
 
 Finally, we let $C'_n= C_n\cup D_1$ and $C'_{n+1}= D_2$ so that $\summation(C'_n)= \sum_{i=1}^{N'_{n}} a_{n,i} +x_{n}$ and $\summation(C'_{n+1})= \sum_{i=1}^{N'_{n+1}} a_{n,i} +x_{n+1}$.
 Observing that $\{b_i\mid 1\le i\le N_n\}\cup D_1\cup D_2$ is the union of the first $N'_n$ terms of $\langle a_{n,i}\rangle_{i=1}^\infty  $ and the first $N'_{n+1}$ terms of $\langle a_{n+1,i}\rangle_{i=1}^\infty $, we conclude that $C_1\cup C_2\cup\cdots \cup C_{n-1}\cup C'_n\cup C'_{n+1}$ is a partition of
 $$\bigcup_{j=1}^{n-1} \{ a_{j,i} \mid 1\le i\le N_j \} \cup \{ a_{n,i}\mid 1\le i\le N'_n  \}\cup \{ a_{n+1,i}\mid 1\leq i\leq N'_{n+1}\}$$
 with the desired property.
 This completes the induction step and therefore the proof of the lemma.
\end{proof}

We are finally ready to prove Theorem~\ref{long-path-forests}, restated with the assumption and notation in this section as follows.

\begin{theorem}\label{030922c}
 Let $n\ge 2$. There exists $L\in \mathbb{N}$ such that every path forest in $\mathcal{H}(n,L)$ is well-burnable.
\end{theorem}   

\begin{proof}
 Fix the number of paths $n$.
 Let $L$ be an integer satisfying the property stated in Lemma~\ref{080322a} for $n-1$.
 By Theorem~\ref{110322c}, it suffices to show that there are no infinite $\prec$-chains starting from a path forest in $\mathcal{H}_{\text{def}}(n,L)$.
 Assume to the contrary that $T=T_0 \prec T_1 \prec T_2 \prec \dotsb$ is an infinite chain starting from some path forest $T$ in $\mathcal{H}_{\text{def}}(n,L)$.
 Along the chain, $\prec$-extensions occur infinitely often at some of the components.
 By starting at some path forest further along the chain if necessary, we may assume without loss of generality that in this infinite chain, \mbox{$\prec$-extensions} occur infinitely often at each of the first $k$ components, while no \mbox{$\prec$-extension} occurs at any of the last $n-k$ components.
 We say that these first $k$ components are \emph{nonstationary} and the rest are \emph{stationary} components.

 Suppose $T=(l_1, l_2, \ldots, l_n)$.
 By the choice of $L$, we can choose an $m\in\mathbb{N}$ so that the path forest $(m^2-\sum_{j=2}^{n} l_j,  l_{2}, l_{3}, \ldots, l_{n})$ is $m$-burnable.
 Of course, $T$ must then have order less that $m^2$ since it is deficient.
 Consider the path forest $T'= (l'_1, l'_2, \ldots, l'_n)$ of order $m^2$ along the infinite chain.
 Since the last $n-k$ components are stationary, we have $(l'_{k+1}, l'_{k+2}, \ldots, l'_n)= (l_{k+1}, l_{k+2}, \ldots, l_n)$.
 If there is only one nonstationary component, then it follows that $l'_1= m^2-\sum_{j=2}^{n} l_j$ and thus $T'$ is $m$-burnable, which gives a contradiction.
 Hence, there are at least two nonstationary components.
 
 Starting from $T'$, the nonstationary components induce $k$ strictly increasing sequences of odd integers without common terms such that the union of all those terms is the set of all odd integers at least $2m+1$.
 Let $\langle a_{j,i}\rangle_{i=1}^\infty$ be the sequence induced by the $j$th component for $1\leq j\leq k$.
 Set $x_j= l'_j$ for each $1\le j\le k-1$ and $x_k= l'_k - \left(m^2-\sum_{j=k+1}^n l'_j\right)$ so that $\sum_{j=1}^k x_j = 0$.
 By Lemma~\ref{070322b}, for some indices $N_1, N_2, \ldots, N_k$, the set $\bigcup_{j=1}^k \{ a_{j,i} \mid 1\le i\le N_j \}$ has a partition $C_1\cup C_2\cup \cdots \cup C_k$ with the property that $\summation(C_j)=  \sum_{i=1}^{N_j} a_{j,i} +x_j$ for every $1\le j\le k$.
 Furthermore, we can choose $N_1, N_2, \dotsc, N_k$ appropriately so that
 $\bigcup_{j=1}^k \{ a_{j,i} \mid 1\leq i\leq N_j \}$ is a set consisting of consecutive odd integers, say up to $2M-1$ for some $M\in\mathbb{N}$.

 Now we consider the larger path forest $T''= (l''_1, l''_2, \ldots, l''_n)$ of order $M^2$ further along the infinite chain.
 We claim that $T''$ is $M$-burnable, which gives us the contradiction we need to complete the proof of the theorem.
 To see this, observe first that by the choice of $m$ and that $(l''_{k+1}, l''_{k+2}, \ldots, l''_n)= (l_{k+1}, l_{k+2}, \ldots, l_n)$, we can place the last $m$ burning sources at the last $n-k+1$ paths of $T''$ in such a way that the last $n-k$ paths are completely burned at the end of the process, while exactly $m^2-\sum_{j=k+1}^n l_j$ vertices of the $k$th path are burned.
 (We could always choose $N_k$ large enough so that the $k$th path of $T''$ has more than $m^2-\sum_{j=k+1}^n l_j$ vertices.)
 
 Note that by the constructions of our sequences, $l''_j=\sum_{i=1}^{N_j} a_{j,i}+l'_j$ for every $1\leq j \leq k$.
 So for each $1\le j\le k-1$, we have $l''_j= \summation(C_j)$, and thus the $j$th path of $T''$ can be burned using the burning sources corresponding to the odd integers in $C_j$.
 Finally, the remaining $l''_k - (m^2-\sum_{j=k+1}^n l_j) = \summation(C_k)$ unburned vertices in the $k$th path of $T''$ can be burned using the burning sources corresponding to the odd integers in $C_k$.
 Therefore, $T''$ is $M$-burnable as claimed.
\end{proof}

\section{Remarks}\label{conclusion}

We have shown in Theorem~\ref{double_m>n} and Theorem~\ref{double_m<n} that the burning number conjecture holds for double spiders.
Moreover, double spiders satisfy the stronger Conjecture~\ref{Tree Conjecture}.
While it will be interesting to see how our work on spiders and double spiders can help in making progress towards the burning number conjecture, we believe the immediate future work is to verify Conjecture~\ref{Tree Conjecture} for the larger class of trees with at most two vertices having degrees greater than two.
This family of trees includes paths, spiders, and double spiders, but we are yet to consider the more general such trees - the union of two spiders, together with a path connecting their respective maximum degree vertices.

\begin{question}
 Suppose $m>n$. Consider a tree $T$ with $n$  leaves of order $m^2+n-2$ such that $T$ has exactly two vertices of degrees at least three.
 Must it be that $T$ is $m$-burnable?
\end{question}

On path forests, our main result in this work shows that every path forest $T$ with a sufficiently long shortest path is well-burnable.
In view of this, we introduce the following definition to study bounds on $L$ in Theorem~\ref{long-path-forests}.

% Comment: While our arguments on path forests are mainly on path orders, but I suppose it is fine to stick with lengths here.

\begin{definition}
 For $n\geq 2$, define $L_n$ to be the least integer with the following property: if $T$ is a path forest with $n$ paths such that its shortest path has order at least $L_n$, then $T$ is well-burnable.
\end{definition}

We know that $L_2 = 3$.
With careful analysis and a little help from a computer, we are also able to determine that $L_3 = 18$ and $L_4 = 26$.
But unfortunately, we do not see how these analyses and arguments can be generalized.

\begin{question}
 What are the values of $L_n$ for $n\ge 5$?
\end{question}

\section{Acknowledgment}

The second author acknowledges the support for this research by the Research University Grant No.~1001/PMATHS/8011129 of Universiti Sains Malaysia.

%\bibliographystyle{abbrv}
%\bibliography{burning}

\begin{thebibliography}{10}
	
	\bibitem{bessy2018bounds}
	S.~Bessy, A.~Bonato, J.~Janssen, D.~Rautenbach, and E.~Roshanbin.
	\newblock Bounds on the burning number.
	\newblock {\em Discrete Appl. Math.}, 235:16--22, 2018.
	
	\bibitem{MR4233796}
	A.~Bonato.
	\newblock A survey of graph burning.
	\newblock {\em Contrib. Discrete Math.}, 16(1):185--197, 2021.
	
	\bibitem{bonato2014burning}
	A.~Bonato, J.~Janssen, and E.~Roshanbin.
	\newblock Burning a graph as a model of social contagion.
	\newblock In {\em Algorithms and Models for the Web Graph}, volume 8882 of {\em
		Lecture Notes in Comput. Sci.}, pages 13--22. Springer, Cham, 2014.
	
	\bibitem{bonato2016how}
	A.~Bonato, J.~Janssen, and E.~Roshanbin.
	\newblock How to burn a graph.
	\newblock {\em Internet Math.}, 12(1-2):85--100, 2016.
	
	\bibitem{bonato2019bounds}
	A.~Bonato and T.~Lidbetter.
	\newblock Bounds on the burning numbers of spiders and path-forests.
	\newblock {\em Theoret. Comput. Sci.}, 794:12--19, 2019.
	
	\bibitem{das2018burning}
	S.~Das, S.~Ranjan~Dev, A.~Sadhukhan, U.~k. Sahoo, and S.~Sen.
	\newblock Burning spiders.
	\newblock In {\em Algorithms and discrete applied mathematics}, volume 10743 of
	{\em Lecture Notes in Comput. Sci.}, pages 155--163. Springer, Cham, 2018.
	
	\bibitem{hiller2019burning}
	M.~Hiller, A.~M. C.~A. Koster, and E.~Triesch.
	\newblock On the burning number of {$p$}-caterpillars.
	\newblock In {\em Graphs and combinatorial optimization: from theory to
		applications---{CTW}2020 proceedings}, volume~5 of {\em AIRO Springer Ser.},
	pages 145--156. Springer, Cham, 2021.
	
	\bibitem{land2016upper}
	M.~R. Land and L.~Lu.
	\newblock An upper bound on the burning number of graphs.
	\newblock In {\em Algorithms and models for the web graph}, volume 10088 of
	{\em Lecture Notes in Comput. Sci.}, pages 1--8. Springer, Cham, 2016.
	
	\bibitem{liu2020burning}
	H.~Liu, X.~Hu, and X.~Hu.
	\newblock Burning number of caterpillars.
	\newblock {\em Discrete Appl. Math.}, 284:332--340, 2020.
	
	\bibitem{mitsche2017burning}
	D.~Mitsche, P.~Pra\l~at, and E.~Roshanbin.
	\newblock Burning graphs: a probabilistic perspective.
	\newblock {\em Graphs Combin.}, 33(2):449--471, 2017.
	
	\bibitem{mitsche2018burning}
	D.~Mitsche, P.~Pra\l~at, and E.~Roshanbin.
	\newblock Burning number of graph products.
	\newblock {\em Theoret. Comput. Sci.}, 746:124--135, 2018.
	
	\bibitem{roshanbin2016burning}
	E.~Roshanbin.
	\newblock {\em Burning a graph as a model of social contagion}.
	\newblock 2016.
	\newblock Thesis (Ph.D.)--Dalhousie University.
	
	\bibitem{tan2020graph}
	T.~S. Tan and W.~C. Teh.
	\newblock Graph burning: tight bounds on the burning numbers of path forests
	and spiders.
	\newblock {\em Appl. Math. Comput.}, 385:125447, 9, 2020.
	
\end{thebibliography}
\end{document}